\documentclass[12pt]{article} 

\usepackage{amsmath}
\usepackage{amsthm}
\usepackage{amsfonts}
\usepackage{mathrsfs}
\usepackage{stmaryrd}
\usepackage{setspace}
\usepackage{fullpage}
\usepackage{amssymb}
\usepackage{breqn}
\usepackage{enumitem}
\usepackage{bbold} 
\usepackage{authblk}
\usepackage{comment}
\usepackage{hyperref}
\usepackage{pgf,tikz}
\usepackage{graphicx}
\usepackage{xcolor}
\usetikzlibrary{decorations.pathreplacing,arrows}
\tikzstyle{vertex}=[circle,draw=black,fill=black,inner sep=0,minimum size=3pt,text=white,font=\footnotesize]

\newtheorem*{rep@theorem}{\rep@title}
\newcommand{\newreptheorem}[2]{%
\newenvironment{rep#1}[1]{%
 \def\rep@title{#2 \ref{##1}}%
 \begin{rep@theorem}}%
 {\end{rep@theorem}}}
\makeatother

\theoremstyle{plain}

\newtheorem*{proposition*}{Proposition}

\newtheorem{thm}{Theorem}[section]

\newtheorem{lemma}[thm]{Lemma}
\newtheorem{conj}[thm]{Conjecture}

\newtheorem{prop}[thm]{Proposition}

\newtheorem{corollary}[thm]{Corollary}

\newtheorem{claim}[thm]{Claim}
\newtheorem{obs}[thm]{Observation}
\newtheorem*{theorem*}{Theorem}
\newtheorem*{prop*}{Proposition}

\newcommand\ex{\ensuremath{\mathrm{ex}}}
\newcommand\rex{\ensuremath{\mathrm{rex}}}
\newcommand\regex{\ensuremath{\mathrm{regex}}}

\newcommand\cG{{\mathcal G}}

\newcommand\cN{{\mathcal N}}

\def\lf{\left\lfloor}   
\def\rf{\right\rfloor}

\newcommand{\ignore}[1]{}

\title{Generalized regular Tur\'an numbers}
\author{D\'aniel Gerbner\footnote{Alfr\'ed R\'enyi Institute of Mathematics, Hun-Ren, Budapest, Hungary. E-mail: \texttt{gerbner@renyi.hu}}, Hilal Hama Karim\footnote{Department of Computer Science and Information Theory, Faculty of Electrical Engineering and Informatics, Budapest University of Technology and Economics, Műegyetem rkp. 3., H-1111 Budapest, Hungary. E-mail: \texttt{hilal.hamakarim@edu.bme.hu}}}

\date{}

\begin{document}

\maketitle

\begin{abstract}
    We combine two generalizations of ordinary Tur\'an problems. Given graphs $H$ and $F$ and a positive integer $n$, we study $\mathrm{rex}(n,H,F)$, which is the largest number of copies of $H$ in $F$-free regular $n$-vertex graphs. We present similar and different behaviours of this function with $\mathrm{ex}(n,H,F)$, which is the maximum number of of copies of $H$ in $n$-vertex $F$-free graphs. Also, we determine the exact value of $\mathrm{rex}(n,K_3,P_k)$, where $P_k$ is a path on $k$ vertices.
\end{abstract}

\section{Introduction}

One of the fundamental theorems in extremal graph theory is due to Tur\'an \cite{turan}. It states that among $n$-vertex $K_{k+1}$-free graphs, the most edges are contained in the complete $k$-partite graph with each part of order $\lfloor n/k\rfloor$ or $\lceil n/k\rceil$. This graph is called the \textit{Tur\'an graph} and is denoted by $T(n,k)$. More generally, given a graph $F$, we denote by $\ex(n,F)$ the largest number of edges in an $n$-vertex $F$-free graph. This topic has attracted a lot of attention, see \cite{fursim} for a survey.

A natural generalization is the study of the largest number of copies of another graph $H$ instead of the number of edges in $n$-vertex $F$-free graphs. Let $\cN(H,G)$ denote the number of not necessarily induced copies of $H$ in $G$, and let $\ex(n,H,F)=\max\{\cN(H,G): \text{$G$ is an $n$-vertex $F$-free graph}\}$. After several sporadic results, 
the systematic study of these so-called generalized Tur\'an problems was initiated by Alon and Shikhelman \cite{alon}.

Another natural generalization was recently considered in \cite{cvk,catu,gptv,Dani1,tati}. Here we study $\rex(n,F)$, which is the largest number of edges in an $n$-vertex \emph{regular} $F$-free graph.

In this paper we combine the above generalizations. Let $\rex(n,H,F):=\max \{\cN(H,G): \text{ $G$ is an $F$-free regular $n$-vertex graph}\}$. Our goal is to show some examples where $\rex(n,H,F)$ behaves similarly to $\ex(n,H,F)$ and also show some examples where they differ significantly.

\smallskip

Given a graph $H$, a blow-up of $H$ is a graph obtained by replacing each vertex of $H$ by an independent set of vertices and each edge $uv \in E(H)$ by a complete bipartite graph between the independent sets replacing $u$ and $v$. A blow-up is balanced, denoted by $H(m)$, if every independent set replacing each vertex is of size $m$, for some $m \in \mathbb{Z}^+$. Alon and Shikhelman \cite{alon} proved the following result.

\begin{thm}[\cite{alon}]\label{blowup}
Let $H$ be a graph on $h$ vertices, then for any graph $F$, we have $\ex(n,H,F)=\Theta(n^{|V(H)|})$ if and only if $F$ is not a subgraph of any blow-up of $H$.
\end{thm}

We extend this theorem to the regular setting.

\begin{thm}\label{teta} For any graph $F$ and $H$, we have that
$\rex(n,H,F)=\Theta(n^{|V(H)|})$ if and only if $F$ is not a subgraph of a blow-up of $H$. 
\end{thm}

Another result of Alon and Shikhelman \cite{alon} is that $\ex(n,K_3,F)=O(n)$ if and only if $F$ is an extended friendship graph. In an extended friendship graph, every cycle is a triangle and there is a vertex $v$ such that every pair of triangles intersect in $v$ (or equivalently, its 2-core is empty or a Friendship graph). We extend this theorem as well to our setting.

\begin{thm}\label{tria}
$\rex(n,K_3,F)=O(n)$ if and only if $F$ is an extended friendship graph.
\end{thm}

Let us turn to problems where adding the regularity changes the situation. It is well-known and easy to see that for any forest $F$, any graph with minimum degree at least $|V(F)|$ contains $F$. This implies that $\rex(n,F)\le (|V(F)|-1)n$. Let $H$ be a connected graph, then the vertices of $H$ have an ordering such that each but the first vertex has a neighbor that is earlier in the ordering. The copies of $H$ in an $F$-free $r$-regular graph can be counted by picking the vertices in the above order. The first vertex can be picked $n$ ways, and then each other vertex can be picked at most $r$ ways among the neighbors of at least one of the vertices picked earlier. This shows that $\rex(n,H,F)=O(n)$, since $r<|V(F)|$. On the other hand, $\ex(n,P_\ell,P_k)=\Theta(n^{\lceil \ell/2\rceil})$ by a theorem of Gy\H ori, Salia, Tompkins and Zamora \cite{GySTZ}.

Another example where the order of magnitude of $\rex(n,H,F)$ is much smaller than that of $\ex(n,H,F)$ is given by even cycles. When $C_{2k}$ is forbidden, the regularity does not have to be constant, but it is $O(n^{1/k})$ by a theorem of Bondy and Simonovits \cite{bs}. Therefore, $\rex(n,C_\ell,C_{2k})=O(n^{1+\frac{\ell-1}{k}})$, while we have $\ex(n,C_\ell,C_{2k})=\Theta(n^{\lfloor \ell/2\rfloor})$ if $3\le \ell\neq 2k$ \cite{gisha}.

Note that we have $\ex(n,C_\ell,C_{2k+1})=\Theta(n^\ell)$ if $\ell$ is even or $\ell>2k+1$, as shown by the blow-up of $C_\ell$. Interestingly, in the remaining case $3<\ell<2k+1$ is odd, we have $\ex(n,C_\ell,C_{2k+1})=\Theta(n^{\lfloor \ell/2\rfloor})$ \cite{gisha}, while the above argument does not give any non-trivial bound. 
It is a natural question to ask whether $\rex(n,C_\ell,C_{2k+1})$ is significantly smaller in this case. We can answer this question in the negative.

\begin{prop}\label{cikli}
   If $3<\ell<2k+1$ is odd, then $\rex(n,C_\ell,C_{2k+1})=\Theta(n^{\lfloor \ell/2\rfloor})$.
\end{prop}

So far we considered only the order of magnitude of $\rex(n,H,F)$. Let us turn to exact and asymptotic results. As shown in \cite{catu,cvk}, for $k\ge3$ we have $\rex(n,K_{k+1})=(1+o(1))|E(T(n,k))|$ (we have $\rex(n,K_{k+1})=|E(T(n,k))|$ if $k$ divides $n$). The exact value of $\rex(n,K_{k+1})$ was determined for all sufficiently large $n$ in \cite{Dani1}.
Let $T^*(n,k)$ denote an arbitrary $n$-vertex $K_{k+1}$-free regular graph with $\rex(n,K_{k+1})$ edges.
Forbidding $K_3$ is very different from forbidding larger cliques in the regular Tur\'an problem. If $n$ is even, then $T(n,2)$ is the regular $n$-vertex triangle-free graph with the most edges. 
If $n$ is odd, then a regular $n$-vertex triangle-free graph with the most edges is obtained by deleting some edges of an $n$-vertex blow-up of $C_5$, as shown in \cite{cvk,catu}.

Given $H$ with $\chi(H)\le k$, there has been a lot of research on whether $\ex(n,H,K_{k+1})=\cN(H,T(n,k))$ for sufficiently large $n$, for example it is the case when $H$ is a complete $l$-partite graph with $3\leq l \leq k$, see e.g. \cite{gerbner,gp,gyps}. There have been two types of counterexamples found (where even $\ex(n,H,K_{k+1})=(1+o(1))\cN(H,T(n,k))$ does not hold). If $H$ is a very unbalanced bipartite graph,
then an unbalanced complete $k$-partite graph may contain more copies of $H$ than the Tur\'an graph. For some graphs $H$, there are $n$-vertex $K_{k+1}$-free graphs that contain more copies of $H$ than any $n$-vertex complete $k$-partite graph. For example, let $H$ be obtained from a path on vertices $v_1,v_2,v_3,v_4,v_5,v_6$ by adding $s$ additional leaves connected to $v_2$ and $s$ additional leaves connected to $v_5$. Then some unbalanced blowup of $C_5$ contains more copies of $H$ than any bipartite graph, see \cite{gyps}. Examples for $k>2$ can be found in \cite{ggyst}. In each of the known constructions, most of the vertices of $H$ would belong to two different classes of $k$-partite graphs, but they can belong to the same class of the blow-up of another graph. Then that class has many vertices.

Both counterexamples are very far from being regular. This suggests that maybe there are no regular counterexamples at all.

\begin{conj}\label{coni}
    Let $k\ge 3$ and $\chi(H)\le k$. Then $\rex(n,H,K_{k+1})=(1+o(1))\cN(H,T(n,k))$. Moreover, if $n$ is sufficiently large and is divisible by $k$, then $\rex(n,H,K_{k+1})=\cN(H,T(n,k))$.
\end{conj}

We prove Conjecture \ref{coni} for complete $k$-partite graphs $H$.

\begin{prop}\label{kpa}
      Let $k\ge 3$ and $H$ be a complete $k$-partite graph. Then $\rex(n,H,K_{k+1})=(1+o(1))\cN(H,T(n,k))$. Moreover, if $n$ is sufficiently large and is divisible by $k$, then $\rex(n,H,K_{k+1})=\cN(H,T(n,k))$.
\end{prop}

We also prove that the moreover part of Conjecture \ref{coni} holds in the case $k=2$. In the case $n$ is odd, the situation is very different, but we can describe the structure of the extremal graph.

\begin{prop}\label{bipa}
    Let $H$ be a bipartite graph. If $n$ is even and sufficiently large, then $\rex(n,H,K_3)=\cN(H,T(n,2))$.
    If $H$ is a tree and $n$ is odd and sufficiently large, then $\rex(n,H,K_3)=\cN(H,G^*)$, where $G^*$ is a regular graph obtained by deleting some edges of an $n$-vertex blow-up of $C_5$.
\end{prop}

Finally, we determine the exact value for $\rex(n,K_3,P_k)$, when $n$ is large enough and $P_k$ is a path on $k$ vertices.
To ease the notation and describe the extremal graphs, we define some graphs first. Let $\cG_{k-1}$ denote the graphs obtained from $K_{k-1}$ by removing the edges of a triangle-free 2-regular subgraph, i.e., the union of vertex-disjoint cycles of length more than 3 such that the total length of the cycles is $k-1$.
In the case $k$ is even, let $G_{k-2}:= K_{k-2}-M$, a clique on $k-2$ vertices in which a perfect matching is removed. Note that each of the above graphs is $(k-4)$-regular and $P_k$-free.  If $k$ is odd, let $G_{k-1}':= K_{k-1}-M$, a clique on $k-1$ vertices in which a perfect matching is removed. Note that $\cN(K_3,G_{k-1})=8 \binom{k/2-1}{3}+3-k/2$ for any graph $G_{k-1}\in \cG_{k-1}$, $\cN(K_3,G_{k-2})=8 \binom{k/2-1}{3}$ and $\cN(K_3,G_{k-1}')=(k-1)(k-3)(k-5)/6=8 \binom{(k-1)/2}{3}$. One can obtain these computations from the fact that the complements of these graphs contain no triangles, the degrees of their vertices and the following formula that is basically proven by Goodman \cite{Gdman}.
$$\mathcal{N}(K_3,G)+ \mathcal{N}(K_3,\overline{G})=\binom{n}{3}-\frac{1}{2}\sum_{v \in V(G)}\mathrm{deg}(v)(n-1-\mathrm{deg}(v)),$$
where $\overline{G}$ denotes the complement of $G$. We denote by $H+F$ the disjoint union of two graphs $H$ and $F$, and by $mF$ we mean $m$ disjoint copies of the graph $F$.

\begin{thm}\label{paths}
    Let $P_k$ be a path on $k$ vertices and $n$ be large enough. Then:

    \begin{enumerate}
        \item If $(k-1)|n$, then $\displaystyle \rex(n,K_3,P_k)=\frac{n}{k-1} \binom{k-1}{3}$, and the unique extremal graph is $\frac{n}{k-1}K_{k-1}$.

        \item Assume that $(k-1) \nmid n$, $k\geq 6$ and either $k-2$ divides $n$ or $k$ is odd.
        Let $n=a(k-2)+b$ with $b<k-2$. Then we have
        $\displaystyle \rex(n,K_3,P_k)=(a-b) \binom{k-2}{3}+ 8 b \binom{\frac{k-1}{2}}{3}$, and the unique extremal graph is $(a-b) (K_{k-2})+ b G'_{k-1}$. 

        \item If $k\geq 6$ is even, and $n$ is neither divisible by $k-1$ nor by $k-2$. Let $n=a(k-3)+b$, with $b<k-3$. Then $$\displaystyle \rex(n,K_3,P_k)=(a-\ell-\lf b/2 \rf) \binom{k-3}{3} + \ell \cN(K_3, G_{k-2})  + \lf b/2 \rf \cN(K_3, G_{k-1}),$$ and the 
        extremal graphs are formed by adding $\lf b/2 \rf$ graphs from $\cG_{k-1}$ to $(a-\ell-\lf b/2 \rf) K_{k-3} + \ell G_{k-2}$,
        where $\ell=0$ if $b$ is even and $\ell=1$ otherwise.

        \item 
        If neither $3$, nor $4$ divides $n$, then $\rex(n,K_3,P_5)=\lfloor n/3\rfloor-1$, and the unique extremal graph is formed by adding a $C_4$ or a $C_5$ to $(\lfloor n/3\rfloor-1)K_3$. 
        In all the cases not listed above, $\rex(n,K_3,P_k)=0$. 
    \end{enumerate}
\end{thm}


\section{Tools}

We will use the following well-known theorem of Erd\H os and Sachs \cite{ersa}.

\begin{thm}[\cite{ersa}]
For every $r$ and $g$, there exists an $r$-regular graph of girth at least $g$.
\end{thm}

In fact we rely on the following simple corollaries of the above theorem.

\begin{lemma}\label{lem1}
   \textbf{(i)} For any $r$ and $k$, if $n$ is sufficiently large and $nr$ is even, then there is an $n$-vertex $r$-regular graph with girth at least $k$.

   \textbf{(ii)} For any $r$, $k$ and $i$, if $n$ is sufficiently large and $nr-i$ is even, then there is an $n$-vertex graph with girth at least $k$ that contains $i$ vertices of degree $r-1$ and each other vertex has degree $r$. Moreover, we can have that the vertices of degree $r-1$ are at distance at least $k-1$.
\end{lemma}

\begin{proof} Let us start by proving \textbf{(i)}.
    We know such a graph $G_1$ exists on $m$ vertices for some $m$. If $r$ is even, we take $r/2$ vertex-disjoint copies of $G_1$ and remove an edge from each. We add a new vertex and connect it to the endpoints of the removed edges. The resulting graph $G_2$ satisfies the desired properties on $\frac{rm}{2}+1$ vertices. For each $n\ge rm^2$, we can write $n$ as $a(\frac{rm}{2}+1)+bm$, thus we can create an $n$-vertex graph by taking vertex-disjoint copies of $G_1$ and $G_2$.

    If $r$ is odd, since $nr$ is even, we must have $n$ is even. In this case, we take $r$ vertex-disjoint copies of $G_1$ and remove an edge from each. We add two new vertices $u,v$ and connect $u$ to one of the endpoints of each removed edge and $v$ to the other endpoint. The resulting graph $G_2'$ satisfies the desired properties on $rm+2$ vertices. For each even $n\ge r^2m$, we can write $n$ as $a(rm+2)+bm$, thus we can create an $n$-vertex graph by taking vertex-disjoint copies of $G_1$ and $G_2'$.

We continue with the proof of \textbf{(ii)}.
    If $i$ is even, we take a graph guaranteed by \textbf{(i)} and remove $i/2$ independent edges such that the endpoints of these edges are at distance at least $k-1$. If $n$ is sufficiently large, we can greedily find such edges. Indeed, we take an edge $u_1v_1$, then at most $2r-2$ other vertices are adjacent to $u$ or $v$, and at most $2(r-1)^j$ vertices are at distance $j$ from $u$ or $v$. Altogether there are at most $2(r-1)^k$ vertices at distance at most $k-1$ from $u$ or $v$. We take a vertex $u_2$ different from those at most $2(r-1)^k$ vertices and an arbitrary neighbor $v_2$ of $u_2$. Repeating this, we can find $i/2$ edges if we can pick a vertex $u_{i/2}$ that is not among the $i-2$ vertices picked earlier and the at most $(i-2)(r-1)^k$ vertices at distance at most $k-1$ from the vertices picked earlier. In other words, we can pick the desired edges if $n>i-2+(i-2)(r-1)^k$. Note that the distance of $u_i$ and $v_i$ is at least $k-1$ after removing the edge $u_iv_i$ because of the girth condition.
    
    
    If $i$ is odd, observe that both $n$ and $r$ are odd. Let $G_3$ be an $(r-1)$-regular $m$-vertex graph for some odd $m$, with girth at least $k$. Let $G_4$ be an $r$-regular graph on $m'$ vertices for some $m'$ sufficiently large with girth at least $k$. Note that $G_3$ and $G_4$ exist by \textbf{(i)}. We take $(m-i)/2$ copies of $G_4$ and remove an edge from each. This way we obtain $m-i$ vertices of degree $r-1$, we connect each of them to a different vertex of $G_3$. The resulting graph has exactly $i$ vertices of degree $r-1$ and each other vertex has degree $r$.
    
    In each of the above constructions, we removed an edge $uv$ from some copy of a graph of girth at least $k$, then we added some edges incident to $u$ and $v$ and outside vertices. After removing $uv$, the distance of $u$ and $v$ becomes at least $k-1$, thus this way we do not create cycles of length less than $k$.
    
    The $i$ vertices of degree $r-1$ (in the copy of $G_3$) can be chosen to be at distance at least $k-1$, by a similar reasoning as in the case $i$ is even.
    
\end{proof}

\begin{corollary}\label{corol}
      For any sequence $(a_n)$ of positive integers with $a_n=\omega(1)$, we can take for every $n$ an $n$-vertex graph $G_n$ that satisfies the assumptions of Lemma \ref{lem1} with $r\le a_n$ and $r=\omega(1)$.
\end{corollary}

The following observation is a simple corollary of Hall's theorem, which will be used in our proofs.

\begin{obs}\label{obsi}
    For every $k\le n$, there exists a $k$-regular bipartite graph with both parts of order $n$.
\end{obs}



A theorem of Andr\'asfai, Erd\H os, and S\'os \cite{aes}, states that a non-bipartite triangle-free graph on $n$ vertices contains a vertex of degree at most $2n/5$. Using this, we prove a stability result on $\rex(2n+1,K_3)$, which may be interesting on its own. When we talk about $V_{i+j}$ in the statement or the proof, then $+$ is meant modulo 5.

\begin{lemma}\label{lem2}
Let $G$ be a $d$-regular $n$-vertex triangle-free graph with $n$ odd. Let $d\ge 2n/5-o(n)$.
Then $V(G)$ contains disjoint sets $V_1,\dots V_5$ such that $|V_i|=n/5-o(n)$ and from $V_i$ there is no edge to $V_i$, $V_{i+2}$ and $V_{i+3}$, and $n/5-o(n)$ edges go to $V_{i+1}$ and $V_{i+4}$. In particular $G$ is obtained by deleting some edges of an $n$-vertex blow-up of $C_5$.
\end{lemma}

\begin{proof} Observe that $G$ cannot be bipartite, thus $d\le 2n/5$ by the result of Andr\'asfai, Erd\H os, and S\'os \cite{aes}.
    Let $C_{2k+1}$ be a shortest odd cycle in $G$ and $C$ be a copy of $C_{2k+1}$. Then every vertex outside $C$ is adjacent to at most two vertices of $C$. This implies that there are at most $2(n-2k-1)$ edges between $C$ and the other vertices. On the other hand, there are at least $(2k+1)d-(2k+1)\ge 2(2k+1)n/5-o(n)$ edges between $C$ and the other vertices by our assumption on the degrees of the vertices of $C$ (which is $d$). Here we use that there are $2k+1$ edges inside $C$, since it is the shortest odd cycle and so does not have any chords.
    
    
    This shows that $k\le 2$. Since $G$ is triangle-free, we have $k=2$. 
    Furthermore, $n-o(n)$ vertices outside $C$ have two neighbors in $C$, otherwise there are at most $2(n-2k-1)-\Omega(n)<(2k+1)d-(2k+1)$ edges between $C$ and the other vertices. Let $v_1,\dots,v_5$ be the vertices of $C$ in the cyclic order. Observe that no vertex can be adjacent to both $v_i$ and $v_{i+1}$, thus $n-o(n)$ vertices are each, for some $i$, adjacent to $v_i$ and $v_{i+2}$. We place those vertices to $V_{i+1}$. Let $U=V(G)\setminus (V_1\cup V_2\cup V_3\cup V_4\cup V_5)$, then $|U|=o(n)$. 


    Let $u\in V_i$. As $u$ has a common neighbor with every vertex of $V_i$, $V_{i+2}$ and $V_{i+3}$, there are no neighbors of $u$ in $V_i\cup V_{i+2}\cup V_{i+3}$, thus all the neighbors of $u$ are in $V_{i+1}$ and $V_{i+4}$ and $U$. In particular, $|V_{i+1}|+|V_{i+4}|\ge 2n/5-o(n)$. This holds for every non-adjacent pair of classes. If $|V_i|\le n/5-\alpha n$, then $|V_{i+2}|,|V_{i+3}|\ge n/5+\alpha n-o(n)$. Then $|V_{i+1}|+|V_{i+4}|\le 2n/5-\alpha n-o(n)$, thus $\alpha=o(1)$, completing the proof.
\end{proof}

 \section{Proofs}

Let us prove Theorem \ref{teta}. Recall that it states that $\rex(n,H,F)=\Theta(n^{|V(H)|})$ if and only if $F$ is not a subgraph of a blow-up of $H$.

\begin{proof}[Proof of Theorem \ref{teta}]
 If $F$ is a subgraph of a blow-up of $H$, then $\rex(n,H,F)\le \ex(n,H,F)=o(n^{|V(H)|})$, where we use Theorem \ref{blowup}. Assume now that $F$ is not a subgraph of any blow-up of $H$. If $H$ is the empty graph, the statement follows. Observe that otherwise, since each bipartite graph is a subgraph of the (sufficiently large) blow-up of a single edge, we also have that $F$ has chromatic number at least 3.
 We can also assume that there are no isolated vertices in $F$.

In the analogous statement for $\ex(n,H,F)$, this is the trivial direction, as the blow-up $H(m)$ with $m=\lfloor n/|V(H)|\rfloor$ is $F$-free and contains $\Omega(n^{|V(H)|})$ copies of $H$. However, we have two problems here: the first is that $H(m)$ is not regular if $H$ is not regular, and the second is that we may need to add some vertices of degree 0 to obtain an $n$-vertex graph. 

Let $\Delta$ be the largest degree in $H$. Let $H_i$ be a graph with girth more than $3|V(F)|$ that has a set $S$ of $i$ vertices of degree $2\Delta$ and all the other vertices of degree $2\Delta+1$, such that the vertices of $S$ are of distance at least $|V(F)|$. Such a graph exists by Lemma \ref{lem1} where the number of vertices is large enough compared to $\Delta$ and $|V(F)|$, 
but constant compared to $n$.

For each $1\leq i \leq \Delta$ and each vertex $v$ of $H$ of degree $i$, we take a copy of $H_{2\Delta+1-i}$ and join $v$ to the $2\Delta+1-i$ vertices of this copy of degree $2\Delta$. This way we obtain a $(2\Delta+1)$-regular graph $H'$ on constant many vertices.

\begin{claim}
The blow-up $H'(m)$ is $F$-free for any $m$. 
\end{claim}

\begin{proof}[Proof of Claim]
We can assume that $m$ is large enough compared to $|V(F)|$. Let us assume that there is a copy of $F$ in $H'(m)$, that we will denote with $F^*$. Let $H^*$ denote an arbitrary copy of an $H_i$ in $H'$, for some $i$. Let $F_0$ denote a connected component of the intersection of $F^*$ with the blow-up of $H^*$. Observe that $F_0$ is bipartite, since any cycle inside $F_0$ has length at most $|V(F)|$ and any cycle inside $H^*$ has length at least $3|V(F)|$. Also observe that $F_0$ contains at most one vertex of $S$. Indeed, otherwise $F_0$ would contain a path between two vertices of $S$ inside $H^*$, but such a path contains more than $|V(F)|$ vertices, a contradiction.
Let $u$ be the vertex of $H$ that is joined to $v$ in $H'$ and $u'$ be an arbitrary neighbor of $u$ in $H$.

Now we can delete $F_0$ and embed it to the complete bipartite graph between the blow-ups of $u$ and $u'$, using only vertices that were not in $F^*$. This can be done since $F_0$ is bipartite and $m$ is large enough. We repeat this for every subgraph of $F$ outside $H(m)$. At the end, we obtain a copy of $F$ in $H(m)$, a contradiction.
\end{proof}

Let us return to the proof of the theorem. We are done if $|V(H')|$ divides $n$, as we can pick $m$ to be $n/|V(H')|$. 
To prove the theorem for every $n$, we do the following. 
Let $H''$ denote the vertex-disjoint union of $H'(2m)$ and $C_{2|V(F)|+1}((2\Delta+1)m)$. Note that $H''$ is $F$-free, since every subgraph of $C_{2|V(F)|+1}((2\Delta+1)m)$ on at most $|V(F)|$ vertices is bipartite. If there is a copy of $F$ in $H''$, then the components that are in $C_{2|V(F)|+1}((2\Delta+1)m)$ are bipartite, and hence could be easily replaced by copies in $H'(2m)$ (we can find such copies in the blow-up of any edge). This way we find a copy of $F$ in $H'(2m)$, a contradiction.

Clearly, $H''$ is $2(2\Delta+1)m$-regular for any $m$, and the number of vertices have the same parity as $m$. Let us pick the largest $m$ such that $n-|V(H'')|$ is even. Observe that $n-|V(H'')|$ is a constant. Now we modify the $C_{2|V(F)|+1}((2\Delta+1)m)$ subgraph. Note that this is similar to the way the odd cycles were modified in \cite{catu}.

Let $A_1,\dots,A_{2|V(F)|+1}$ be the blown-up parts of the cycle in this order.
We take a pair of neighboring parts, say $A_i$ and $A_{i+1}$, and add $b=(n-|V(H'')|)/2$ vertices to each of $A_i$ and $A_{i+1}$. We add them in such a way that we still have a complete bipartite graph between any pair of consecutive blown-up parts $A_j,A_{j+1}$, i.e., we connect the new vertices of $A_i$ to each vertex of $A_{i-1}$ and $A_{i+1}$, and connect the new vertices of $A_{i+1}$ to each vertex of $A_i$ and $A_{i+2}$. Then we remove the edges of a spanning bipartite graph $B$ between $A_{i-1}$ and $A_i$ such that each vertex of $A_{i-1}$ has degree $b$ and each vertex of $A_i$ has degree $\left\lfloor\frac{(2\Delta+1)mb}{(2\Delta+1)m+b}\right\rfloor$ or $\left\lceil\frac{(2\Delta+1)mb}{(2\Delta+1)m+b}\right\rceil$ in $B$. 
This can be done 
the following way. We cover $A_i$ by $|A_{i-1}|$ sets of size $b$ (each of them are the $b$ neighbors of a vertex in $A_{i-1}$), such that each vertex of $A_{i}$ is covered by as equal as possible number of such sets. Then we take a bijection between the vertices of $A_{i-1}$ and these sets, and delete the edges between the vertices of $A_{i-1}$ and the vertices of the sets they are mapped to.

We remove the edges of a copy of $B$ between $A_{i+1}$ and $A_{i+2}$ as well such that the vertices of degree $b$ are in $A_{i+2}$. 

At this point the vertices outside $A_i$ and $A_{i+1}$ have degree $2(2\Delta+1)m$. The part $A_i$ consists of a set $A_i'$ of vertices of degree $2(2\Delta+1)m+b-\left\lfloor\frac{(2\Delta+1)mb}{(2\Delta+1)m+b}\right\rfloor$ and a set $A_i''$ of vertices of degree $2(2\Delta+1)m+b-\left\lceil\frac{(2\Delta+1)mb}{(2\Delta+1)m+b}\right\rceil$. Similarly, $A_{i+1}$ consists of a set $A_{i+1}'$ of vertices of degree $2(2\Delta+1)m+b-\left\lfloor\frac{(2\Delta+1)mb}{(2\Delta+1)m+b}\right\rfloor$ and a set $A_{i+1}''$ of vertices of degree $2(2\Delta+1)m+b-\left\lceil\frac{(2\Delta+1)mb}{(2\Delta+1)m+b}\right\rceil$. Observe that by the analogous construction, we have that $|A_i'|=|A_{i+1}'|$. We pick a perfect matching $M'$ between $A_i'$ and $A_{i+1}'$, and extend it to a perfect matching $M$ between $A_i$ and $A_{i+1}$. We delete the edges of $M$.

Then the resulting graph between $A_i$ and $A_{i+1}$ is $2(2\Delta+1)m+b-1$-regular, thus we can delete matchings between $A_i$ and $A_{i+1}$ till we obtain a $\left((2\Delta+1)m+\left\lfloor\frac{(2\Delta+1)mb}{(2\Delta+1)m+b}\right\rfloor\right)$-regular graph between $A_i$ and $A_{i+1}$. After that, we add the edges of $M$ that are not in $M'$. Observe that vertices of $A_i'$ have $\left((2\Delta+1)m+\left\lfloor\frac{(2\Delta+1)mb}{(2\Delta+1)m+b}\right\rfloor\right)$ neighbors in $A_{i+1}$ and $(2\Delta+1)m-\left\lfloor\frac{(2\Delta+1)mb}{(2\Delta+1)m+b}\right\rfloor$ neighbors in $A_{i-1}$. Vertices of $A_i''$ have one more neighbor in $A_{i+1}$ and one less neighbor in $A_{i-1}$. The same holds for vertices in $A_{i+1}$.
Let $G$ denote the resulting $n$-vertex graph. Then $G$ is $2(2\Delta+1)m$-regular and contains at least $m^{|V(H)|}=\Theta(n^{|V(H)|})$ copies of $H$, completing the proof. 
\end{proof}

Let us continue with the proof of Theorem \ref{tria}. Recall that it states
that
$\rex(n,K_3,F)=O(n)$ if and only if $F$ is an extended friendship graph.

\begin{proof}[Proof of Theorem \ref{tria}]
If $F$ is an extended friendship graph, then $\rex(n,K_3,F)\le \ex(n,K_3,F)=O(n)$.

Assume that $F$ is not an extended friendship graph. Then it either contains 
two vertex-disjoint triangles or a longer cycle $C_k$ with $k\ge 4$. In the first case, we take the $K_3$-free graph on $n-1$ vertices with regularity $r=\Omega(n)$ due to Caro and Tuza \cite{catu}. In particular, it contains an induced copy $K$ of $K_{r/2,r/2}$. We remove a perfect matching from $K$, and add a new vertex $v$, connected to the vertices of $K$. The resulting graph is $r$-regular, and contains $r^2/4$ triangles that all contain $v$, completing our proof.

Let us assume now that $F$ contains $C_k$. Let $n$ be sufficiently large. We pick $r=\omega(1)$ such that $r$ is small enough to have a $2r$-regular graph of girth more than $3k^2$ on $m$ vertices whenever $m\ge 11n/36$ and such that $n>r^{3|V(F)|}$. We take an $r$-regular $m$-vertex graph $G_0$ of girth at least $3k^2$ where $m=\lfloor n/4\rfloor$, using Corollary \ref{corol}. We consider $G_0^k$ as an auxiliary graph. Recall that the $k$th power $G_0^k$ of a graph $G_0$ is obtained by joining vertices of distance at most $k$.

It is easy to see that $G_0^k$ is $r'$-reg, where $r'=r+r(r-1)+\dots+r(r-1)^{k-1}$. We take a proper $r'+1$-edge-coloring of $G_0^k$. Since $G_0$ is a subgraph of $G_0^k$, we obtain a proper edge-coloring of $G_0$. For each color $i$, we partition the edges of color $i$ to some number of $r$-sets and a set of order at most $r$. For each such set, we add a new vertex and connect it to the at most $2r$ vertices that are incident to those at most $r$ edges. This way we obtain $G_1$. 

The vertices of $G_0$ have degree $r$ in $G_0$, thus they are incident to edges of $r$ colors, hence their degree is $2r$ in $G_1$. The newly added vertices have degree $2r$, except $r'+1$ vertices, that are connected to the endpoints of less than $r$ edges. Let us assume that the sum of degrees in $G_1$ is $2r|V(G_1)|-\ell$. Note that $\ell$ is even since each vertex has an even degree and at most $2r(r'+1)$. There are $rm$ edges from $V(G_0)$ to $V(G_1)\setminus V(G_0)$ and at least $2r(|V(G_1)\setminus V(G_0)|-r'-1)$ edges from $V(G_1)\setminus V(G_0)$ to $V(G_0)$, thus $|V(G_1)\setminus V(G_0)|\le \frac{m}{2}+r'+1$, hence $V(G_1)\le 3n/8+O(1)\le 4n/9$.


Now we make $G_1$ regular. We take a copy of a $2r$-regular graph $G_0'$ and remove $\ell/2$ independent edges the following way. First we delete an arbitrary edge $u_1v_1$. The number of vertices at distance at most $|V(F)|+1$ from $u_1$ is at most $r''=2r+2r(2r-1)+\dots+2r(2r-1)^{|V(F)|+1}$. Then we pick a vertex $u_2$ that is at distance at least $|V(F)|+1$ from $u_1$ and a neighbor $v_2$ of $u_2$. Then we delete the edge $u_2v_2$ and these four vertices are at distance at least $|V(F)|$ from each other. We repeat this, always picking vertices $u_i$ that are at distance at least $|V(F)|+2$ from each of $u_1,\dots,u_{i-1}$. This is doable if $\ell r''<m$, which holds by our assumption on $r$.

The resulting graph $Q_1$ is of girth more than $k$ with a set $S$ of $\ell$ vertices of degree $2r-1$ and all the other vertices of degree $2r$. The vertices of $S$ are at distance at least $|V(F)|$.
We join each vertex $v$ of $G_1$ to $2r-d(v)$ vertices of degree $2r-1$ in this new graph. The resulting graph $G_2$ is $2r$-regular on at most $25n/36$ vertices. 

Finally, we add a $2r$-regular graph of girth more than $k$ on $n-|V(G_2)|$ vertices. This exists by the choice of $r$.
\end{proof}

Let us continue with the proof of Proposition \ref{cikli}. Recall that it states that if $3<\ell<2k+1$ is odd, then $\rex(n,C_\ell,C_{2k+1})=\Theta(n^{\lfloor \ell/2\rfloor})$.

\begin{proof}[Proof of Proposition \ref{cikli}]
The upper bound is shown by $\ex(n,C_\ell,C_{2k+1})=\Theta(n^{\lfloor \ell/2\rfloor})$ in \cite{gisha}.

Let us turn to the lower bound. We start with an unbalanced blow-up of $C_\ell$, where we blow up $(\ell-1)/2$ independent vertices to $m$-sets, and keep the other vertices (note that this construction shows the analogous bound for $\ex(n,C_\ell, C_{2k+1})$, but it is far from regular). Let $H$ denote this graph, then the largest degree is $2m$ in $H$. We take two vertex-disjoint copies of $H$, denoted by $2H$. We add sets $A_1,A_2,\dots,A_{2k}$ of new vertices of order $2m-2$. We take two blown-up parts $A$ and $A'$ of order $m$ of $2H$ arbitrarily. We add sets $A_1,A_2,\dots,A_{2k}$ of new vertices of order $2m-2$. We take all the possible edges between $A$ and $A_1$, then an $m$-regular graph between $A_i$ and $A_{i+1}$ for each $i\le 2k-1$ (this exists because of Observation \ref{obsi}), and then take all the possible edges between $A_{2k}$ and $A'$. It is easy to see that each vertex of $A$, $A'$ and each $A_i$ has degree $2m$ and no $C_{2k+1}$ is created this way. We repeat this by taking $4k(m-1)$ new vertices as long as there are at least two blown up classes of order $m$ in $2H$.

We are left with several vertices of degree $2m$ and exactly two adjacent vertices $u,v$ and $u',v'$ of degree $m+1$ in both copies of $H$ (they are the adjacent vertices of the original $\ell$-cycles that were not blown up). We take sets $B$, $B'$ of order $m-1$ and $B_1,\dots, B_{2k}$ of order $2m-1$. We take all the edges between $u$ and $B$ and between $B$ and $B_1$. Then we take an $(m+1)$-regular graph between $B_{2i+1}$ and $B_{2i+2}$, for each $0\leq i\leq k-1$, and an $(m-1)$-regular graph between $B_{2i}$ and $B_{2i+1}$, for each $1\leq i \leq k-1$. Finally, we take all the edges between $B_{2k}$ and $B'$ and between $B'$ and $v$. We do the same to deal with $u'$ and $v'$ in the other copy of $H$.





It is left to add $n-|V(H')|$ vertices without ruining these properties. Observe that we added at most $16\ell km$ vertices to $H$. We pick $m$ to be the largest odd number below $\lfloor n/40k\ell\rfloor$, thus $H'$ has at most $n/2$ vertices. If $n-|V(H')|$ is even, we can pick a bipartite $2m$-regular graph on those vertices, completing the proof. If $n-|V(H')|$ is odd, we additionally pick a copy of $C_{2k+3}(m)$, and then pick a bipartite $2m$-regular graph on the remaining vertices, completing the proof.
\end{proof}

Let us continue with the proof of Proposition \ref{kpa}. Recall that it states that if 
$H$ is a complete $k$-partite graph
and $k\ge 3$, then $\rex(n,H,K_{k+1})=(1+o(1))\cN(H,T(n,k))$, and without the error term if $k$ divides $n$. Let $T^*(n,k)$ denote an arbitrary extremal graph for $\rex(n,K_{k+1})$. Recall that if $k|n$, then $T^*(n,k)$ is
$T(n,k)$, otherwise it is obtained from $T(n,k)$ by changing a small number of edges, as shown in \cite{catu,cvk}.

\begin{proof}[Proof of Proposition \ref{kpa}] The lower bound follows from $T^*(n,k)$, since as mentioned above it is very close to the Tur\'an graph. We will prove the upper bound for $k\ge 2$. Note that if $k=2$ the case $n$ is even is equivalent to the even case of Proposition \ref{bipa}.

Let $H=K_{s_1,\dots,s_k}$ with $s_1\le s_2\le\dots \le s_k$.
    For simplicity, we will deal with labeled copies of $H$, clearly the same $n$-vertex $K_{k+1}$-free regular graph maximizes (asymptotically) the number of labeled copies of $H$ as the number of copies. We will show that for $k\ge 2$, any $n$-vertex $K_{k+1}$-free graph contains at most $k!\prod_{i=1}^k\binom{(n/k)}{s_i}$ copies of $H$. Clearly $T(n,k)$ satisfies this with equality if $k$ divides $n$, and $T^*(n,k)$ gives the correct asymptotics, thus this upper bound completes the proof for other values of $n$ as well.

We apply induction on $k$ and on $\sum_{i=1}^k s_i$. The base case $k=2$ follows from Proposition \ref{bipa} if $n$ is even. Moreover, the exact same proof works for the case $n$ is odd, without changing a single word (but gives a bound that is not sharp). The other base case $\sum_{i=1}^k s_i=k$ follows from Zykov's theorem \cite{zykov}, which states that $\ex(n,K_r,K_{k+1})=\cN(H,T(n,k))$.


  Let $G$ be an $n$-vertex $K_{k+1}$-free $r$-regular graph, then $r\le (k-1)n/k$ by Tur\'an's theorem. We consider two cases. Assume first that $s_1=1$ and let $H'$ be the graph we obtain by deleting the first class from $H$. Then we first pick a vertex $v$ of $G$ corresponding to the single vertex in the first class, at most $n$ ways. Then we pick a labeled copy of $H'$ in the neighborhood of $v$, at most $(k-1)!\prod_{i=2}^k\binom{(r/(k-1))}{s_i}$ by the induction on $k$. This way we picked the labeled copies of $H$ at most $n(k-1)!\prod_{i=2}^k\binom{(r/(k-1))}{s_i}\le nk!/k\prod_{i=2}^k\binom{(k-1)n/k)/(k-1)}{s_i}=k!\frac{n}{k}\prod_{i=2}^k\binom{n/k}{s_i}=k!\binom{n/k}{s_1}\prod_{i=2}^k\binom{n/k}{s_i}$.

  Assume now that $s_1>1$ and let $H''$ be the graph we obtain by deleting a $K_k$ from $H$. We first pick an unlabeled copy $K$ of $K_k$, then a labeled copy of $H''$ from the remaining vertices, and then add the labels to the vertices of $K$. By Zykov's theorem, the number of unlabeled copies of $K_k$ is maximized by the Tur\'an graph, thus it is at most $n^k/k^k$. The number of labeled copies of $H''$ is at most $k!\prod_{i=1}^k\binom{(n-k)/k}{s_i-1}$ by induction on $\sum_{i=1}^k s_i$. Afterwards, we add the vertices of $K$ to the vertices of $H''$. Observe that each vertex of $K$ has a copy of $K_{k-1}$ in their neighborhood in $H''$. As $G$ is $K_{k+1}$-free, the vertices of a $K_{k-1}$ cannot be adjacent to two adjacent vertices. This implies that each copy of $K_{k-1}$ has at most one common neighbor in $K$. Therefore, each vertex of $K$ can belong to at most one of the classes of $H$. This means that the number of labels the vertices of $K$ can receive is at most $\prod_{i=1}^k s_i$, hence  
  the number of labeled copies of $H$ is at most $\frac{1}{\prod_{i=1}^k s_i}\frac{n^k}{k^k} k!\prod_{i=1}^k\binom{(n-k)/k}{s_i-1}= k!\prod_{i=1}^k\binom{n/k}{s_i}$. This completes the proof.
\end{proof}

Let us continue with Proposition \ref{bipa}. Recall that it extends the above proposition to the case $k=2$ if $n$ is even. If $n$ is odd, then it deals with the case $H$ is a tree and claims that the extremal graph $G^*$ is a regular graph obtained by deleting some edges of a blow-up of $C_5$.


\begin{proof}[Proof of Proposition \ref{bipa}]
Let $n$ be even, consider a component $H_0$ of $H$ and an ordering of the vertices of $H_0$ such that each but the first vertex has an earlier neighbor. Such an ordering  exists by first picking an arbitrary vertex and then each time picking a neighbor of a vertex already picked.

Let $G$ be an $n$-vertex $r$-regular triangle-free graph, then $r\le n/2$. Moreover, either $G$ is bipartite (with both parts of order $n/2$, thus $G$ is contained in $T(n,2)$, completing the proof), or $r\le 2n/5$. Assume that $r\le 2n/5$. There are at most $n$ ways to pick the first vertex and at most $2n/5$ ways to pick each subsequent vertex. In $T(n,2)$ there are $n$ and $(1+o(1))n/2$ ways to do this. The copies of $H_0$ may be counted multiple times, but the number of times is a fixed constant $c$ depending only on the automorphisms of $H$ and not the host graph. Therefore, $G$ contains at most $2\left(\frac{2n}{5}\right)^{|V(H_0)|}/c$, while $T(n,2)$ contains $(2+o(1))\left(\frac{n}{2}\right)^{|V(H_0)|}/c$ copies of $H_0$. Then we pick the other components of $H$. Similarly, there are more ways to pick each component in $T(n,2)$ than in $G$ if $n$ is large enough. Therefore, $T(n,2)$ contains more copies of $H$ for sufficiently large $n$, completing the proof.

    Assume now that $n$ is odd and $H$ is a forest. Let $G$ be an $n$-vertex $r$-regular triangle-free graph, then $r\le 2n/5$. Moreover, by Lemma \ref{lem2}, either $G$ is obtained by deleting some edges of an $n$-vertex blow-up of $C_5$, or $r\le 2n/5-\varepsilon n$ for some $\varepsilon>0$. In the first case, we are done. In the second case, we can proceed similarly to the argument in the case where $n$ is even. $G$ contains 
    at most $2\left(\frac{2n}{5}-\varepsilon n)\right)^{|V(H_0)|}/c$ copies of $H_0$, while $G^*$ contains $2\left(\frac{2n}{5}\right)^{|V(H_0)|}/c$ copies of $H_0$. The same holds for other components, thus $G^*$ contains more copies of $H$ for sufficiently large $n$, completing the proof in this case.
\end{proof}

Before the proof of Theorem \ref{paths}, let us mention some results that will be used. By a theorem of Erd\H{o}s and Gallai \cite{erga}, if a connected graph has at least $k$ vertices and minimum degree $\lf k/2 \rf$, then it contains a $P_k$. Gerbner, Patk\'os, Tuza and Vizer \cite{Dani1} gave the exact value of $\regex(n,T)$ for any tree $T$ and large $n$. Note that $\regex(n,F)= \max \{d \ : \ \text{there is an} \ n\text{-vertex}, \ d\text{-regular}, \ F\text{-free graph} \ G \}$

\begin{thm}[\cite{Dani1}]\label{regexT}
 Let $T$ be a tree on $t$ vertices and $n>n_0(T)$. Then
 
 $\operatorname{regex}(n, T)= \begin{cases}t-2 & \text { if } t-1 \text { divides } n \text { or } T \text { is a star and } t \text { or } n \text { is even, } \\ t-3 & \text { if the above does not hold, and either } t \text { is odd, or } t-2 \text { divides } n, \\ & \text { or } T \text { is a star or } T \text { is an almost-star and } n \text { is even. } \\ t-4 & \text { otherwise. }\end{cases}$

\end{thm}

A tree is  an \textit{almost-star}
if in its proper 2-coloring, one of the classes consists of at most two vertices (thus a path on at least 6 vertices is not an almost-star).

Now we are ready to present the proof of Theorem \ref{paths}, which determines $\rex(n,K_3,P_k)$ if $k\ge 7$ and $n$ is sufficiently large.



\begin{proof}[Proof of Theorem \ref{paths}]

Let $G$ be an $n$-vertex $r$-regular $P_k$-free graph containing the maximum number of triangles. Observe that each vertex of $G$ is in at most $\binom{r}{2}$ triangles, and hence $\cN(K_3,G)\leq \frac{n}{3} \binom{r}{2}$, with equality only when each vertex is in a clique $K_{r+1}$. Note that in each of the constructions described in the introduction, the $r$-regular $n$-vertex graph contains $n/(r+1)-O(1)$ copies of $K_{r+1}$, hence contains $\frac{n}{3} \binom{r}{2}-O(1)$ triangles. This implies that, for large $n$, a graph with smaller regularity cannot contain more triangles than our construction. Therefore, $r$ is at least the regularity of the claimed unique construction. Using Theorem \ref{regexT}, in each of the cases we know that $r$ is at most the regularity of the claimed unique construction, and it is left to show that no other $r$-regular graph can contain at least as many triangles as our construction, in each case.

As $G$ is $P_k$-free, we have  $r \leq k-2$ by Theorem \ref{regexT}. 
Therefore, $\cN(K_3,G) \leq \frac{n}{3}\binom{k-2}{2}=\frac{n}{k-1}\binom{k-1}{3}$, with equality only when $(k-1)|n$ and $G$ is $n/k-1$ disjoint copies of $K_{k-1}$, proving the first case. If $n$ is not divisible by $(k-1)$ and $r=k-2$, then by Theorem \ref{regexT} we have that $P_k$ is a star, i.e., $k\le 3$.
Consequently, $r < k-2$.

In the second case, we have $k$ is odd or $k-2$ divides $n$ and we can assume $r=k-3$. We can write $n$ as $a(k-2)+b$, where $0 \leq b \leq k-3$. Since $k\geq 6$, we have $k-3\ge k/2$, and hence, by the result of Erd\H{o}s and Gallai, each component has at most $k-1$ vertices. Thus, each component is either a $K_{k-2}$ or $G_{k-1}'$, for these are the only $(k-3)$-regular graphs on at most $k-1$ vertices. This means $G=xK_{k-2}+yG_{k-1}'$, which gives $x(k-2)+y(k-1)= n =a(k-2)+b$, implying 
$$y=a'(k-2)+b, \ \text{where} \ a'=a-x-y.$$ If $y>b$, then $a'\geq 1$, and hence, $y\geq (k-2)+b$, then we can replace $k-2$ copies of $G_{k-1}'$ by $k-1$ copies of $K_{k-2}$, increasing the number of triangles as $\cN(K_3, K_{k-2}) > \cN(K_3, G_{k-1}')$, contradicting the choice of $G$. Also, if $y<b$, then as $k-2>0$, we must have $a'<0$, which implies $y<0$, a contradiction. Therefore, $y=b$ and $x=a-b$, which proves the second case.

Note that if $k$ is even, then $K_{k-1}$ does not contain a perfect matching, and hence, if $(k-2) \nmid n$, then $r < k-3$, leading to the third case.

Assume now $r=k-4$ and $k\geq 6$ is even. We can write $n$ as $a(k-3)+b$, where $0 \leq b \leq k-4$. First recall that each graph in $\cG_{k-1}$ contains the same number of triangles.
By the same reasoning of the previous case, we may assume that $G$ consists of $x K_{k-3}+y G_{k-2}$ and $z$ copies of graphs from $\cG_{k-1}$.
This gives $x(k-3)+y(k-2)+z(k-1)=n=a(k-3)+b$, implying 
$$2z\geq a'(k-3)-y+b, \ \text{where} \ a'=a-x-y-z.$$
Note that $\cN(K_3,K_{k-3})+\cN(K_3,G_{k-1}) > 2 \cN(K_3, G_{k-2})$, for any graph $G_{k-1}\in \cG_{k-1}$, and hence, whenever there are two copies of $G_{k-2}$ in $G$, we can replace them by a copy of $K_{k-3}$ and a copy of $G_{k-1}$, increasing the number of triangles. Therefore, we have that $y$ is either $0$ or $1$.
If $z>\lf b/2 \rf$, then we have $a'\geq 1$, which means $z \geq  (k-3)/2 -y/2 +  b/2$. If $b$ is even we may assume $y=0$, and hence, in both cases of $b$ being odd or even, we still have $\displaystyle z\geq \lceil (k-3)/2 \rceil + \lf b/2 \rf$. We can then replace $\lceil (k-3)/2\rceil$ copies of $G_{k-1}$ by $\lf (k-1)/2 \rf$ copies of $K_{k-3}$ and a copy of $G_{k-2}$, increasing the number of triangles, which contradicts the extremality of $G$. Again, due to compatibility of the number of vertices, $z$ cannot be less than $\lf b/2 \rf$, proving the third case.

Finally, if $k=5$, the only connected regular $P_k$-free graph with regularity at least 3 is $K_4$. In the case of regularity 2, the only connected 2-regular graph that contains a triangle is $K_3$. If 3 does not divide $n$, we need to add at least one longer cycle to a graph consisting of vertex-disjoint triangles. If $k=4$, the only connected regular $P_k$-free graph that contains a triangle is $K_3$, thus if $n$ not divisible by 3, then $\rex(n,P_4)=0$. In the remaining cases $k\le 3$, the triangle contains $P_k$, thus $\rex(n,P_k)=0$, completing the proof.
\end{proof}

\bigskip

\textbf{Funding}: Research supported by the National Research, Development and Innovation Office - NKFIH under the grants FK 132060 and KKP-133819.

\end{document}